\documentclass[12pt]{article}

\usepackage{amsmath,amssymb,amsthm}
\usepackage[top=1in,bottom=1in,left=1in,right=1in]{geometry}
\usepackage{graphicx}

\newtheorem{theorem}{Theorem}[section]
\theoremstyle{definition}
\newtheorem{example}{Example}[section]
\newtheorem{remark}{Remark}[section]

\newcommand{\real}{\mathbb{R}}
\newcommand{\parder}[2]{\frac{\partial {#1}}{\partial {#2}}}

\begin{document}

\begin{center}
\textbf{\Large Patchy Solution of a Francis\textendash Byrnes\textendash Isidori\\ Partial  
Differential  Equation\footnote{Research supported in part by AFOSR and NSF.} }\\
Cesar O. Aguilar\footnote{National Research Council Postdoctoral Fellow, Department of Applied Mathematics, Naval Postgraduate School, Monterey, CA 93943, \texttt{coaguila@nps.edu}} and Arthur J. Krener\footnote{Distinguished Visiting Professor, Department of Applied Mathematics, Naval Postgraduate School, Monterey, CA 93943, \texttt{ajkrener@nps.edu}}
\end{center}

\begin{abstract}
The solution to the nonlinear output regulation problem requires one to solve a first order PDE, known as the Francis-Byrnes-Isidori (FBI) equations. In this paper we propose a method to compute approximate solutions to the FBI equations when the zero dynamics of the plant are hyperbolic and the exosystem is two-dimensional.  With our method we are able to produce approximations that converge uniformly to the true solution.  Our method relies on the periodic nature of two-dimensional analytic center manifolds.
\end{abstract}

\baselineskip 1.5em

\section{Introduction}\label{sec:intro}
Consider the control system 
\begin{equation}\label{eqn:full-sys}
\begin{aligned}
\dot x &= f(x,u,w)\\
\dot w &= s(w)\\
y &= h(x,u,w)
\end{aligned}
\end{equation}
where $x\in\real^n$ is the state variable, $u\in\real^m$ is the control variable, $w\in\real^q$ is an exogenous variable, and $y\in\real^p$ is the output variable.  The maps $f:\real^n\times\real^m\times\real^q\rightarrow\real^n$, $s:\real^q\rightarrow\real^q$ and $h:\real^n\times\real^m\times\real^q\rightarrow\real^p$ are all assumed to be sufficiently smooth and satisfy $f(0,0,0)=0$, $s(0)=0$, and $h(0,0,0)=0$.  The variable $w$ represents a disturbance and/or a reference signal, and its dynamics are commonly referred to as the \textit{exosystem}.  The \textit{state feedback regulator problem} \cite{IsiByr90} is to find a state feedback control $u=\alpha(x,w)$, with $\alpha(0,0)=0$, such that the equilibrium $x=0$ of the dynamical system
\[
\dot x = f(x,\alpha(x,0),0)
\]
is exponentially stable, and such that for each sufficiently small initial condition $(x_0,w_0)$  the solution of \eqref{eqn:full-sys} with $u=\alpha(x,w)$ satisfies
\[
\lim_{t\rightarrow\infty} y(t) =0.
\]

A characterization of the state feedback regulator problem for linear systems was given by Francis \cite{Fran77} and later generalized to nonlinear systems by Isidori and Byrnes \cite{IsiByr90}.  As shown in \cite{IsiByr90}, the solvability of the regulator problem can be reduced to the solvability of a system of partial differential equations (PDEs), which in the linear case reduce to the Sylvester type equation obtained by Francis.  For this reason, we refer to these equations as the Francis\textendash Byrnes\textendash Isidori (FBI) PDEs.  For completeness, we state the main result of \cite{IsiByr90} (for the definition of Poission stability used below see Remark~\ref{rem:neutral}).
\begin{theorem}[Isidori\textendash Byrnes]
Assume that the equilibrium $w=0$ of the exosystem is Lyapunov stable and there is a neighborhood of $w=0$ in which every point is Poisson stable.  Assume further that the pair
\begin{equation}\label{eqn:lin-sys}
\left(\parder{f}{x}(0,0,0), \parder{f}{u}(0,0,0)\right)
\end{equation}
is stabilizable.  Then the state feedback regulator problem is solvable if and only if there exists $C^k$ ($k\geq 2$) mappings $\pi:\Omega\rightarrow\real^n$, with $\pi(0)=0$, and $\kappa:\Omega\rightarrow\real^m$, with $\kappa(0)=0$, both defined in a neighborhood $\Omega\subseteq\real^q$ of $w=0$, and satisfying 
\begin{equation}\label{eqn:fbi-eqns}
\begin{aligned}
\parder{\pi}{w}(w) s(w) &= f(\pi(w),\kappa(w),w)\\
0&= h(\pi(w),\kappa(w),w).
\end{aligned}
\end{equation}
\end{theorem} 
Given a solution pair $(\pi,\kappa)$ to the FBI equations \eqref{eqn:fbi-eqns}, a state feedback solving the regulator problem is given by
\[
\alpha(x,w)= \kappa(w) + K(x-\pi(w))
\]
where $K\in\real^m\times\real^n$ is any feedback matrix rendering the pair \eqref{eqn:lin-sys} asymptotically stable.  In general, solutions to the FBI equations, being singular quasilinear PDEs with constraints, may not exist.  However, for a class of control-affine systems, it is shown in \cite{IsiByr90} that the solvability of the FBI equations is a property of the zero dynamics of \eqref{eqn:full-sys}.  Roughly speaking, if the zero dynamics of \eqref{eqn:full-sys} has a hyperbolic equilibrium at the origin then a solution to the FBI equations exists by the center manifold theorem \cite{Carr81}.  It is known, however, that center manifolds suffer from a number of subtle properties associated with uniqueness and differentiability \cite{Sijbrand85}.  Despite these difficulties, a $C^\infty$ dynamical system possess a $C^k$ center manifold for each $k\geq 1$, and moreover, it is possible to obtain approximate solutions of arbitrarily high-order via Taylor series \cite{Carr81}.  In this respect,  Huang and Rugh \cite{HuRu92} and Krener \cite{Kr92} provide a method to compute approximate solutions to the FBI equations via Taylor polynomials, yielding approximate output regulation.  A shortcoming of this approach is that the domain on which the series approximation yields satisfactory results is not guaranteed to enlarge significantly by computing higher order approximations.  This can be a serious drawback as the number of monomials in $q$ variables of degree $d$ is $\binom{q+d-1}{d}$, a number growing rapidly in $d$.  Moreover, polynomial approximations to $(\pi,\kappa)$ can lead to destabilizing effects when the order of the approximation increases.

In this paper, we present a method to compute solutions to the FBI equations for the class of real analytic SISO control-affine systems 
\begin{equation}\label{eqn:cas}
\begin{aligned}
\dot x &=f(x)+g(x)u\\
\dot w &=s(w)\\
y &=h(x)+p(w)
\end{aligned}
\end{equation}
where $f:\real^n\rightarrow\real^n$, $g:\real^n\rightarrow\real^n$, $s:\real^q\rightarrow\real^q$, $h:\real^n\rightarrow\real$, and $p:\real^q\rightarrow\real$ are real analytic about the origin.  Furthermore, we will restrict our considerations to two-dimensional exosystems, i.e., $q=2$, whose linear part contains non-zero eigenvalues.  Our method is based on the existence and uniqueness results for two-dimensional analytic center manifolds in \cite{Aulbach85} and on the high-order patchy method of Navasca and Krener \cite{NavKr07}.  A key strength of our approximation method, which relies on the periodic nature of the solution to the FBI equations for \eqref{eqn:cas}, is the reduction of the computational effort inherent in a direct Taylor polynomial approximation.  

The organization of this paper is as follows. In Section~\ref{sec:ana-sol-fbi} we briefly summarize the key insight provided in \cite{IsiByr90} on how the solvability of the FBI equations can be reduced to the problem of solving a center manifold equation provided the zero dynamics of \eqref{eqn:cas} are hyperbolic.  With this simplification, we show how the standard stability assumptions on the exosystem lead to a direct application of the results in \cite{Aulbach85} to deduce uniqueness of solutions to the FBI equations of \eqref{eqn:cas}.  In Section~\ref{sec:algorithm}, we describe our method to compute a high-order piecewise smooth approximation to the solution of the FBI equations and prove that a sequence of approximations generated by our method converges uniformly to the true solution.  Finally, in Section~\ref{sec:examples} we illustrate our method on examples and then make some concluding remarks.

\begin{remark}\label{rem:neutral}
Henceforth, it will be implicitly assumed that the exosystem has an equilibrium $w=0$ that is Lyapunov stable and that there is a neighborhood of $w=0$ in which every point is Poisson stable.  We will refer to this type of stability as \textit{neutral stability}.  By Poisson stability we mean the following.  An initial condition $x_0$ of the dynamical system $\dot x=f(x)$ is Poisson stable if the flow $\Phi^f_t(x_0)$ of the vector field $f$ is defined for all $t\in\real$ and for each neighbourhood $U$ of $x_0$ and for each real number $T>0$, there exists a time $t_1>T$ such that $\Phi^f_{t_1}(x_0)\in U$ and a time $t_2<-T$ such that $\Phi^f_{t_2}(x_0)\in U$.
\end{remark}

%=====================================================
\section{Real analytic and periodic solutions to the FBI equations}\label{sec:ana-sol-fbi}
As shown in \cite{IsiByr90}, a key simplification in the problem of solving the FBI equations for a system of the form \eqref{eqn:cas} consists in reducing it to the problem of solving a center manifold equation for the zero dynamics of \eqref{eqn:cas}.  Following \cite{IsiByr90} and using the now standard notation in \cite{IsiBook}, assume that the triplet $\{f,g,h\}$ has relative degree $1\leq r<n$ at $x=0$ and let $(z,\xi)$ denote the standard normal coordinates, where $\xi=(h(x),L_fh(x),\ldots,L_f^{r-1}h(x))$ and $z$ is such that $L_gz=0$.  In the $(z,\xi)$ coordinates, \eqref{eqn:cas} takes the form
\begin{equation}\label{eqn:new-sys}
\begin{aligned}
\dot z &= f_0(z,\xi)\\
\dot \xi_1 &= \xi_2, \ldots, \dot \xi_{r-1} = \xi_r\\
\dot \xi_r &= b(z,\xi) + a(z,\xi)u\\
\dot w &= s(w)\\
y &= \xi_1 + p(w).
\end{aligned}
\end{equation}
The zero dynamics of \eqref{eqn:cas} are given by the dynamical system
\begin{equation}\label{eqn:zero-dyn}
\dot z = f_0(z,0). 
\end{equation}
Define functions $\varphi_i:\real^q\rightarrow\real$ by $\varphi_i(w)=-L_s^{i-1}p(w)$, $1\leq i\leq r$, set $\varphi(w)=(\varphi_1(w),\ldots,\varphi_{r}(w))$, and let
\[
u_e(x,w) = -\frac{L^r_fh(x)+L_s^rp(w)}{L_gL_f^{r-1}h(x)}.
\]
Then it is straightforward to verify that if $\phi$ satisfies the PDE
\begin{equation}\label{eqn:reduced}
\parder{\phi}{w}(w) s(w) = f_0(\phi(w),\varphi(w))
\end{equation}
then $\pi(w):=(\phi(w),\varphi(w))$ and $\kappa(w):=u_e(\pi(w),w)$ constitute a solution pair to the FBI equations of \eqref{eqn:new-sys}.  If the origin of \eqref{eqn:zero-dyn} is hyperbolic, then \eqref{eqn:reduced} is the equation that is satisfied by any center manifold $\{(z,w): z=\phi(w)\}$ of the dynamical system
\begin{equation}\label{eqn:zero-exo}
\begin{aligned}
\dot z &= f_0(z,\varphi(w))\\
\dot w &= s(w).
\end{aligned}
\end{equation}
Hence, in the hyperbolic case, the problem of solving the FBI equations associated to the original system \eqref{eqn:cas} is reduced to solving the center manifold equation associated to \eqref{eqn:zero-exo}.  Although this simplification is significant, solutions to center manifolds suffer from subtleties associated with uniqueness and differentiability \cite{Sijbrand85}.  For example, it is known that an analytic dynamical system does not generally posses an analytic center manifold, thereby forcing one to seek a center manifold solution that is only $C^k$ ($k=2,3,\ldots$) and thus not necessarily unique.  As an example, the polynomial dynamical system 
\begin{equation*}
\begin{aligned}
\dot z &= -z + w_1^2+w_2^2\\
\dot w_1 &= -w_2-\tfrac{1}{2}w_1(w_1^2+w_2^2)\\
\dot w_2 &= w_1-\tfrac{1}{2}w_2(w_1^2+w_2^2)
\end{aligned}
\end{equation*} 
which is of the form \eqref{eqn:zero-exo}, has the property that each center manifold has Taylor series
\[
\sum_{i=1}^\infty (i-1)!(w_1^2+w_2^2)^i
\]
which has vanishing radius of convergence.  Despite these difficulties, a special case for which sharp uniqueness and differentiability results exist is for two-dimensional center manifolds, and is given by the following theorem due to Aulbach \cite{Aulbach85}.
\begin{theorem}[Aulbach]\label{thm:aulbach}
Consider the ordinary differential equation
\begin{equation}\label{eqn:aulbach}
\begin{aligned}
\dot z &= Bz + Z(w_1,w_2,z)\\
\dot w_1 &= -w_2 + P(w_1,w_2,z)\\
\dot w_2 &= w_1 + Q(w_1,w_2,z)
\end{aligned}
\end{equation}
where $w_1,w_2\in\real$, $z\in\real^n$, and $P,Q,$ and $Z$ are real analytic functions about the origin and have Taylor series beginning with quadratic terms.  Suppose that the matrix $B$ has no eigenvalues on the imaginary axis.  If the local center manifold dynamics of \eqref{eqn:aulbach} are Lyapunov stable and non-attractive then \eqref{eqn:aulbach} has a uniquely determined local center manifold which is analytic and generated by a family of periodic solutions.
\end{theorem}
Aulbach's result has a direct application to the output regulation problem, as given by the following theorem.
\begin{theorem}\label{eqn:exi-cm}
Suppose that in \eqref{eqn:cas} the exosystem is two-dimensional and $\frac{\partial s}{\partial w}(0)$ has non-zero eigenvalues.  Suppose that $f,g,h$ and $p$ are real analytic mappings about $x=0$ and $w=0$, respectively, and that the triple $\{f,g,h\}$ has a well-defined relative degree $1\leq r<n$ at $x=0$.  If the zero dynamics of \eqref{eqn:cas} are hyperbolic, then there exist unique and real analytic mappings $(\pi,\kappa)$ solving the associated FBI equations of \eqref{eqn:cas}.
\end{theorem}
\noindent\textbf{Proof.}  By assumption and neutral stability of the exosystem, the eigenvalues of the exosystem are non-zero and purely imaginary.  Indeed, if the eigenvalues were not purely imaginary then $w=0$ would necessarily be either a repelling or an attractive equilibrium, contradicting the assumption of neutral stability.  Now since $B:=\parder{f_0}{z}(0,0)$ contains eigenvalues off the imaginary axis, there exists an analytic coordiante change \cite{Perko91} about the origin such that \eqref{eqn:zero-exo} takes the form 
\begin{equation}\label{eqn:exo-zero-coor}
\begin{aligned}
\dot z &= Bz + Z(w_1,w_2,z)\\
\dot w_1 &= -w_2 + P(w_1,w_2)\\
\dot w_2 &= w_1 + Q(w_1,w_2),
\end{aligned}
\end{equation}
where $P,Q$ and $Z$ are analytic at the origin and have Taylor series beginning with quadratic terms.  From \eqref{eqn:exo-zero-coor}, we can observe that the dynamics of any center manifold of \eqref{eqn:exo-zero-coor} are equivalent to the exosystem dynamics, which by assumption are Lyapunov stable and non-attractive.  Aulbach's theorem completes the proof.\hfill$\blacksquare$

\begin{remark}
Theorem~\ref{eqn:exi-cm} actually holds for more general MIMO control-affine systems with $m=p$.  In \cite{Hu03}, it is shown that if the composite control-affine system
\begin{align*}
\dot x &= f(x,w) + \sum_{i=1}^m g_i(x,w) u_i\\
\dot w &= s(w)\\
y &= h(x,w)
\end{align*}
has a well-defined relative degree at $(x,w)=(0,0)$, then the associated FBI equations are solvable if the zero dynamics of the composite system are hyperbolic.  In this case, the FBI equations reduce to a center manifold equation of the form \eqref{eqn:reduced} so that Aulbach's theorem can be applied when the exosystem is two-dimensional.  
\end{remark}

\begin{example}\label{exm:pendulum}
The dynamics of a cart and inverted pendulum system can be written in the form 
\begin{equation}\label{eqn:cart-pend}
\begin{aligned}
\dot x_1 &= x_2\\
\dot x_2 &= u\\
\dot x_3 &= x_4\\
\dot x_4 &= \frac{g}{\ell}\sin(x_3) - \frac{1}{\ell}\cos(x_3)u
\end{aligned}
\end{equation}
where $x_1$ is the position of the cart, $x_3$ is the angle the pendulum makes with the vertical, $g$ is the acceleration due to gravity, $\ell$ is the length of the rod, and $u$ is the control force.  With $h(x)=x_1$, the system has relative degree $r=2$ at $x=0$, and therefore $(\xi_1,\xi_2)=\xi(x)=(h(x),L_fh(x))=(x_1,x_2)$.  With $(z_1,z_2) = z(x)=(x_3, x_4+\tfrac{x_2}{\ell}\cos(x_3))$,  the zero dynamics are given by
\begin{align*}
\dot z_1 &= z_2\\
\dot z_2 &= \frac{g}{\ell}\sin(z_1),
\end{align*}
whose linearization has eigenvalues $\pm\sqrt{\tfrac{g}{\ell}}$.  Hence, with system output $y=x_1+p(w)$ ($p$ real analytic) and a two-dimensional real analytic exosystem whose linearization has non-zero eigenvalues, there exists a unique and real analytic solution to the associated FBI equations of the cart and inverted pendulum system \eqref{eqn:cart-pend}.\hfill$\square$
\end{example}

%===============================================
\section{Computation of the center manifold}\label{sec:algorithm} 
In this section we outline a method to compute the solution to the FBI equations in the case of two-dimensional exosystem and real analytic data.  As described in the previous section, for the nonlinear control systems in consideration, the solvability of the FBI equations can be reduced to solving a center manifold equation for a dynamical system of the form
\begin{equation}\label{eqn:gen-cm}
\begin{aligned}
\dot z &= Bz + \bar{Z}(w_1,w_2,z)\\
\dot w_1 &= -w_2 + P(w_1,w_2)\\
\dot w_2 &= w_1 + Q(w_1,w_2),
\end{aligned}
\end{equation}
where $w=(w_1,w_2)\in\real^2$, $z\in\real^n$, $\bar{Z}$, $P$, and $Q$ are real analytic mappings, and the eigenvalues of $B$ have non-zero real parts.  We will therefore limit our considerations to solving the center manifold equation for \eqref{eqn:gen-cm}.  It will be assumed that the $w$-dynamics have $w=0$ as a Lyapunov stable and non-attractive equilibrium.  By Theorem~\ref{thm:aulbach}, there exists a unique analytic mapping $\phi(w_1,w_2)$, defined locally about $w=0$, solving the center manifold PDE associated to \eqref{eqn:gen-cm}.

Our method is best described on the representation of \eqref{eqn:gen-cm} in polar coordinates.  Hence, we apply the tranformation $(w_1,w_2,z)=(r\cos\theta,r\sin\theta,z)$ to \eqref{eqn:gen-cm} yielding a system of the form
\begin{equation}\label{eqn:polar-sys}
\begin{aligned}
\dot r &= r\hat{R}(\theta,r)\\
\dot\theta &= 1+\hat\Theta(\theta,r)\\
\dot z &= Bz + \hat{Z}(\theta,r,z)
\end{aligned}
\end{equation}
where $\hat{R},\hat{\Theta},\hat{Z}$ are analytic functions converging for each $\theta\in[0,2\pi]$ and $|r|\leq a$, $\|z\|\leq a$, where $a>0$ is a positive constant.  Define $\hat{f}(\theta,r,z)=Bz + \hat{Z}(\theta,r,z)$.  The center manifold PDE for \eqref{eqn:polar-sys} is 
\begin{equation}\label{eqn:pde-polar}
\hat{f}(\theta,r,\psi(\theta,r)) = \parder{\psi}{\theta}[1+\hat\Theta(\theta,r)] + \parder{\psi}{r} r\hat{R}(\theta,r)
\end{equation}
for the unknown analytic mapping $\psi(\theta,r)$ ($=\phi(r\cos\theta,r\sin\theta)$).  The mapping $\psi$ has a power series representation
\[
\psi(\theta,r)=\sum_{i=1}^\infty e_i(\theta)r^i
\]
converging in a cylinder of the form $\theta\in[0,2\pi]$, $|r|\leq \epsilon$, and with $2\pi$-periodic coefficients $e_i(\theta)$ \cite{Aulbach85}.  By eliminating the time variable $t$, \eqref{eqn:polar-sys} can be reduced to
\begin{subequations}\label{eqn:polar-red}
\begin{align}
\frac{dr}{d\theta} &= rR(\theta,r)\label{eqn:polar-red-r}\\[2ex]
\frac{dz}{d\theta} &= Bz + Z(\theta,r,z)\label{eqn:polar-red-z}.
\end{align}
\end{subequations}
Define $f(\theta,r,z)=Bz+Z(\theta,r,z)$.  From \eqref{eqn:pde-polar} it follows that
\begin{equation}\label{eqn:pde-polar-red}
f(\theta,r,\psi(\theta,r)) = \parder{\psi}{\theta} + \parder{\psi}{r}\frac{dr}{d\theta}.
\end{equation}

We now give a brief sketch of our method.  Let $r(\theta)$ be a solution to \eqref{eqn:polar-red-r} and define the mapping
\[
\Psi(\theta,\sigma)=\psi(\theta,r(\theta)+\sigma)
\]
for $\theta\in[0,2\pi]$ and $|\sigma|$ small.  We note that, with $r=r(\theta)$ substituted into the RHS of \eqref{eqn:polar-red-z}, the curve $\Psi(\theta,0)=\psi(\theta,r(\theta))$ is the solution to \eqref{eqn:polar-red-z} with initial condition $z(0)=\psi(0,r(0))$.  For $|\sigma|$ sufficiently small, we have a power series representation
\begin{equation}\label{eqn:psi-pwr-ser}
\Psi(\theta,\sigma)=\Psi(\theta,0) + \sum_{i=1}^\infty \parder{^i\Psi}{\sigma^i}(\theta,0)\frac{\sigma^i}{i!}
\end{equation}
converging for all $\theta\in[0,2\pi]$ and having $2\pi$-periodic coefficients $\parder{^i\Psi}{\sigma^i}(\theta,0)$.  In fact, it is easy to see that
\begin{equation}\label{eqn:coeffs-Psi-per}
\parder{^i\Psi}{\sigma^i}(\theta,0) = \parder{^i\psi}{r^i}(\theta,r(\theta)).
\end{equation}
By construction, the mapping $\Psi$ is a perturbation of $\psi(\theta,r(\theta))$ in the radial direction, the amount of perturbation given by the parameter $\sigma$.  Our method is based on computing the Taylor series approximation 
\[
\Psi^N(\theta,\sigma)= \Psi(\theta,0) + \sum_{i=1}^N \parder{^i\Psi}{\sigma^i}(\theta,0)\frac{\sigma^i}{i!}
\]
and using it to build the center manifold along $r(\theta)$ in the radial direction.  Having followed $\Psi^N$ along a small annular region, say of the form
\[
\{(\theta,r):0\leq\theta\leq 2\pi, r(\theta)\leq r < r(\theta)+\epsilon\},
\]
we compute a new radial curve $\theta\mapsto \tilde r(\theta)$ with initial condition $\tilde r(0)=r(0)+\epsilon$, compute the new corresponding Taylor series approximation $\tilde\Psi^N$, and then continue building the center manifold by following $\tilde\Psi^N$ along the annular region
\[
\{(\theta,r):0\leq\theta\leq 2\pi, \tilde r(\theta)\leq r < \tilde r(\theta)+\tilde\epsilon\}.
\]
This process is repeated and the annular regions, along with the corresponding approximations, are patched together to form a piecewise smooth approximation to the true solution $\psi$.

\begin{remark}
To compute the Taylor series approximations $\Psi^N$ it is necessary to compute the $\theta$-dependent coefficients appearing in \eqref{eqn:psi-pwr-ser}, which can be done in the following way.  From the definition of $\Psi$, a direct computations gives 
\[
\parder{\Psi}{\theta} = \parder{\psi}{\theta}(\theta,r(\theta)+\sigma) + \parder{\psi}{r}(\theta,r(\theta)+\sigma)\frac{dr}{d\theta}
\]
which when combined with \eqref{eqn:pde-polar-red} yields
\begin{equation}\label{eqn:psi-pde}
\parder{\Psi}{\theta} = f(\theta,r(\theta)+\sigma,\Psi(\theta,\sigma)).
\end{equation}
Using \eqref{eqn:psi-pde}, we can now write a down linear inhomogeneous ODE for the coefficient $\parder{^i\Psi}{\sigma^i}(\theta,0)$.  Indeed, differentiating \eqref{eqn:psi-pde} with respect to $\sigma$, and interchanging the order of differentiation, yields
\[
\parder{}{\theta}\left(\parder{\Psi}{\sigma}(\theta,\sigma)\right) = \parder{f}{z}(\theta,r(\theta)+\sigma,\Psi(\theta,\sigma))\parder{\Psi}{\sigma}(\theta,\sigma) + \parder{f}{r}(\theta,r(\theta)+\sigma,\Psi(\theta,\sigma))
\]
and therefore
\[
\parder{}{\theta}\left(\parder{\Psi}{\sigma}(\theta,0)\right) = A(\theta)\parder{\Psi}{\sigma}(\theta,0) + \parder{f}{r}(\theta,r(\theta),\Psi(\theta,0))
\]
where the matrix $A(\theta)=\parder{f}{z}(\theta,r(\theta),\Psi(\theta,0))$.  In general, it can be verified by induction that 
\begin{equation}\label{eqn:lin-ode-psi}
\parder{}{\theta}\left(\parder{^i\Psi}{\sigma^i}(\theta,0)\right) = A(\theta)\parder{^i\Psi}{\sigma^i}(\theta,0) + F_i\left(\theta,\Psi(\theta,0),\parder{\Psi}{\sigma}(\theta,0),\ldots,\parder{^{i-1}\Psi}{\sigma^{i-1}}(\theta,0)\right)
\end{equation}
for some mappings $F_i$, $i\geq 2$.
\end{remark}

With the previous constructions in mind, we are now ready to describe an algorithm for computing the solution $\psi$ to the center manifold equation \eqref{eqn:pde-polar-red}.  
\begin{enumerate}
\item Let $N\geq 1$ be a fixed positive integer and let 
\[
\psi^N_0(\theta,r)= \sum_{i=1}^N e_i(\theta)r^i,
\]
that is, $\psi^N_0$ is simply the $N$th order Taylor approximation of $\psi$ in $r$.  To compute $\psi^N_0$, one can use the method in \cite{HuRu92} to generate a $N$th order Taylor polynomial approximation of $\phi(w_1,w_2)$, say $\phi^N(w_1,w_2)$, and then simply $\psi^N_0(\theta,r)=\phi^N(r\cos\theta,r\sin\theta)$.  Set $\Psi_0=\psi$ and set $r_{-1}(\theta)=0$ for $\theta\in\real$.  The initial approximation $\psi^N_0$ will be accepted in an annular region of the form
\[
\{(\theta,r):0\leq\theta\leq 2\pi, 0\leq r<r_0(\theta)\}
\]
where $r_0(\theta)$ is the solution to \eqref{eqn:polar-red-r} with some prescribed initial condition $r_0(0)=\epsilon_0>0$.  To compute accurate numerical solutions to $r_0$, we solve a BVP using \eqref{eqn:polar-red-r} with boundary conditions $r(0)=r(2\pi)=\epsilon_0$ and constant initial guess $\epsilon_0$ on $[0,2\pi]$.  

\item Define $\Psi_1(\theta,\sigma)=\psi(\theta,r_0(\theta)+\sigma)$.  From \eqref{eqn:psi-pwr-ser}, $\Psi_1$ can be approximated by the truncated series 
\[
\Psi^N_{1}(\theta,\sigma)=\Psi_1(\theta,0) + \sum_{i=1}^N \parder{^i\Psi_1}{\sigma^i}(\theta,0)\frac{\sigma^i}{i!}
\]
for $|\sigma|$ small.  To obtain accurate numerical solutions to the coefficients $\parder{^i\Psi_1}{\sigma^i}(\theta,0)$, we solve BVPs using the ODEs \eqref{eqn:lin-ode-psi} with boundary conditions $\parder{^i\Psi_1}{\sigma^i}(0,0)=\parder{^i\Psi_1}{\sigma^i}(2\pi,0)$ and initial guesses
\[
\parder{^i\Psi_1}{\sigma^i}(\theta,0) \approx \parder{^i\psi^N_0}{r^i}(\theta,r_0(\theta)).
\]
Similarly, to compute $\Psi_1(\theta,0)=\psi(\theta,r_0(\theta))$ we solve a BVP using \eqref{eqn:polar-red-z} with boundary conditions $z(0)=z(2\pi)$ and initial guess 
\[
\Psi_1(\theta,0) \approx \psi^N_0(\theta,r_0(\theta)).
\]
Having computed $\Psi_1(\theta,0),\parder{\Psi_1}{\sigma}(\theta,0),\ldots,\parder{^N\Psi_1}{\sigma^N}(\theta,0)$, we obtain an approximation $\psi_1(\theta,r)$ to $\psi(\theta,r)$ defined by
\[
\psi_1(\theta,r)=\Psi^N_{1}(\theta,r-r_0(\theta))
\]
which is accepted in the region
\begin{equation}\label{eqn:psi-dom-prv}
\{(\theta,r):0\leq\theta\leq 2\pi, \, r_0(\theta)\leq r < r_0(\theta)+\epsilon_1\}
\end{equation}
for some desired $\epsilon_1>0$.  In this way, we have extended our original approximation $\psi_0$ of $\psi$ to the domain \eqref{eqn:psi-dom-prv}.  Our running approximation of $\psi$ is given by
\[
\psi(\theta,r)\approx \begin{cases} \psi_0(\theta,r), & 0\leq r < r_0(\theta),\\[2ex]
\psi_1(\theta,r), & r_0(\theta)\leq r \leq r_0(\theta)+\epsilon_1
\end{cases}
\]
for $\theta\in[0,2\pi]$.

\item We now proceed to augment to our running approximation a mapping $\psi_2$, that will be defined on an an annular region surrounding the domain of $\psi_1$, in the following way.  We first compute the solution $r_1(\theta)$ to \eqref{eqn:polar-red-r} with initial condition $r_1(0)=r_0(0)+\epsilon_1$.  As in Step 2, this is done by solving a BVP using \eqref{eqn:polar-red-r} with boundary conditions $r(0)=r(2\pi)=r_0(0)+\epsilon_1$ and taking the curve $r_0(\cdot)+\epsilon_1$ as an initial guess to $r_1$.  Here we note that, to avoid overlapping domains of definition between $\psi_1$ and $\psi_2$, the domain \eqref{eqn:psi-dom-prv} of $\psi_1$ is redefined to be
\[
\{(\theta,r):0\leq\theta\leq 2\pi,r_0(\theta)\leq r < r_1(\theta)\}. 
\]

\item We now repeat Step 3 with $r_1(\theta)$ and build an approximation to $\Psi_2(\theta,\sigma)=\psi(\theta,r_1(\theta)+\sigma)$ of the form
\[
\Psi^N_{2}(\theta,\sigma)=\Psi_2(\theta,0) +  
\sum_{i=1}^N \parder{^i\Psi_2}{\sigma^i}(\theta,0)\frac{\sigma^i}{i!},
\]
for $0\leq\sigma\leq\epsilon_2$ and $\epsilon_2>0$ sufficiently small.  The coefficients $\parder{^i\Psi_2}{\sigma^i}(\theta,0)$ are computed by solving BVPs using the ODEs \eqref{eqn:lin-ode-psi} with boundary conditions $\parder{^i\Psi_2}{\sigma^i}(0,0)=\parder{^i\Psi_2}{\sigma^i}(2\pi,0)$ and initial guesses
\[
\parder{^i\Psi_2}{\sigma^i}(\theta,0)\approx \parder{^i\Psi_1}{\sigma^i}(\theta,0).
\]
That is, we use the previously computed coefficients as initial guesses for the current coefficients.  Similarly, to compute $\Psi_2(\theta,0)=\psi(\theta,r_1(\theta))$ we solve a BVP using \eqref{eqn:polar-red-z} with boundary conditions $z(0)=z(2\pi)$ and initial guess 
\[
\Psi_2(\theta,0) \approx \psi_1(\theta,r_1(\theta)).
\]
Having computed $\Psi_2(\theta,0),\parder{\Psi_2}{\sigma}(\theta,0),\ldots,\parder{^N\Psi_2}{\sigma^N}(\theta,0)$, we obtain the approximation
\[
\psi_2(\theta,r)=\Psi^N_{2}(\theta,r-r_1(\theta))
\]
which is accepted in the region
\begin{equation}\label{eqn:psi2-dom-prv}
\{(\theta,r):0\leq\theta\leq 2\pi, \, r_1(\theta)\leq r < r_1(\theta)+\epsilon_2\}.
\end{equation}
In this way, we extend our approximation of $\psi(\theta,r)$ to the annulus \eqref{eqn:psi2-dom-prv}
and our running approximation is
\[
\psi(\theta,r)\approx \begin{cases} \psi_0(\theta,r), & 0\leq r < r_0(\theta), \\[2ex]
\psi_1(\theta,r), & r_0(\theta)\leq r < r_1(\theta), \\[2ex]
\psi_2(\theta,r), & r_1(\theta)\leq r \leq r_1(\theta)+\epsilon_2
\end{cases}
\]
for $\theta\in[0,2\pi]$.

\item  Steps 4-5 can now be iterated.  Indeed, suppose we have computed an approximation $\psi_k(\theta,r)$ to $\psi(\theta,r)$ of the form $\psi_k(\theta,r)=\Psi^N_{k}(\theta,r-r_{k-1}(\theta))$, and defined in the region 
\begin{equation}\label{eqn:psik-prv-dom}
\{(\theta,r):0\leq\theta\leq 2\pi, r_{k-1}(\theta)\leq r < r_{k-1}(\theta)+\epsilon_k\}, 
\end{equation}
where $\Psi_k(\theta,\sigma)=\psi(\theta,r_{k-1}(\theta)+\sigma)$,
\[
\Psi^N_{k}(\theta,\sigma) = \Psi_k(\theta,0) + \sum_{i=1}^N \parder{^i\Psi_k}{\sigma^i}(\theta,0)\frac{\sigma^i}{i!}
\]
and $r_{k-1}(\theta)$ is the solution to \eqref{eqn:polar-red-r} with initial condition $r_{k-1}(0)=\Sigma_{i=0}^{k-1}\epsilon_i$.  To extend our current approximation of $\psi$ beyond the domain of $\psi_k$, we begin by computing the solution $r_k(\theta)$ to \eqref{eqn:polar-red-r} with initial condition $r_k(0)=r_{k-1}(0)+\epsilon_k$.  This is done by solving a BVP using \eqref{eqn:polar-red-r} with boundary conditions $r(0)=r(2\pi)=r_{k-1}(0)+\epsilon_k$ and initial guess $r_{k-1}(\cdot)+\epsilon_k$.  Let now $\Psi_{k+1}(\theta,\sigma)=\psi(\theta,r_k(\theta)+\sigma)$ and form
\[
\Psi^N_{k+1}(\theta,\sigma) = \Psi_{k+1}(\theta,0) + \sum_{i=1}^N \parder{^i\Psi_{k+1}}{\sigma^i}(\theta,0)\frac{\sigma^i}{i!}.
\]
The coefficients $\parder{^i\Psi_{k+1}}{\sigma^i}(\theta,0)$ are computed by solving BVPs using the ODEs \eqref{eqn:lin-ode-psi} with boundary conditions $\parder{^i\Psi_{k+1}}{\sigma^i}(0,0)=\parder{^i\Psi_{k+1}}{\sigma^i}(2\pi,0)$ and initial guesses
\[
\parder{^i\Psi_{k+1}}{\sigma^i}(\theta,0)\approx \parder{^i\Psi_k}{\sigma^i}(\theta,0).
\]
Similarly, to compute $\Psi_{k+1}(\theta,0)$ we solve a BVP using \eqref{eqn:polar-red-z} with boundary conditions $z(0)=z(2\pi)$ and initial guess 
\[
\Psi_{k+1}(\theta,0) \approx \psi_k(\theta,r_k(\theta)).
\]
Having computed $\Psi_{k+1}(\theta,0),\parder{\Psi_{k+1}}{\sigma}(\theta,0),\ldots,\parder{^N\Psi_{k+1}}{\sigma^N}(\theta,0)$, we obtain the approximation
\[
\psi_{k+1}(\theta,r)=\Psi^N_{k+1}(\theta,r-r_k(\theta))
\]
which is accepted in the region
\begin{equation}\label{eqn:psikpp-dom-prv}
\{(\theta,r):0\leq\theta\leq 2\pi, \, r_k(\theta)\leq r < r_k(\theta)+\epsilon_{k+1}\}.
\end{equation}
To avoid overlapping the domains of $\psi_k$ and $\psi_{k+1}$, the domain of definition \eqref{eqn:psik-prv-dom} of $\psi_k$ is redefined to be
\[
\{(\theta,r):0\leq\theta\leq 2\pi, r_{k-1}(\theta)\leq r\leq r_{k}(\theta)\}.
\]

\item After iterating Steps 4-5 a $k\geq 1$ number of times, we obtain the following piecewise smooth approximation to $\psi(\theta,r)$:
\begin{equation}\label{eqn:the-approx}
\psi(\theta,r) \approx \tilde\psi_k(\theta,r):=\begin{cases}
\psi_0(\theta,r), & 0\leq r<r_0(\theta), \\[2ex]
\psi_1(\theta,r), & r_0\leq r<r_1(\theta), \\[2ex]
\quad\vdots & \qquad\qquad\quad\vdots\\
\psi_k(\theta,r), & r_{k-1}\leq r\leq r_k(\theta)
\end{cases}
\end{equation}
for $\theta\in[0,2\pi]$.
\end{enumerate}

Let us make a few remarks about our algorithm.
\begin{remark}
\begin{enumerate}

\item The coefficients $\parder{^i\Psi}{\sigma^i}(\theta,0)$, $i=0,1,\ldots$, are computed by solving BVP problems for mainly two reasons: (1) they are known \textit{a priori} to be periodic, and (2) we have good approximations of them from the previously computed coefficients.  

\item The main computational effort of our method is in computing the coefficients $\parder{^i\Psi}{\sigma^i}(\theta,0)$.  From \eqref{eqn:lin-ode-psi}, we see that the ODE for $\parder{^i\Psi}{\sigma^i}(\theta,0)$ is linear in $\parder{^i\Psi}{\sigma^i}(\theta,0)$ and is polynomial in the previously computed $\Psi(\theta,0),\ldots,\parder{^{i-1}\Psi}{\sigma^{i-1}}(\theta,0)$.  Hence, a way to speed up the computation is to solve for the coefficients $\parder{^i\Psi}{\sigma^i}(\theta,0)$ order-by-order.  This can result in computational savings when the RHS of \eqref{eqn:lin-ode-psi} is complicated to evaluate or when $n\times N$ is large.  

\item When the $w$-dynamics are given by the harmonic oscillator
\begin{align*}
\dot w_1 &= -w_2\\
\dot w_2 &= w_1
\end{align*}
the computation of the radial curves $r(\theta)$ is trivial and are given by the constant curves $r(\theta)=r(0)$.  In this case, there is no need to redefine the outer boundary of the successive approximations $\psi_i$ when going from one annulus to the other.

\item In order for the algorithm to produce a meaningful approximation to $\psi$, the domain on which the approximation \eqref{eqn:the-approx} is defined, namely $\{(\theta,r):0\leq\theta\leq 2\pi, 0\leq r\leq r_{k}(\theta)\}$, must of course be contained in the cylinder $\theta\in[0,2\pi]$, $|r|\leq\epsilon$ on which $\psi$ is defined.  Since $\epsilon$ is not known \textit{a priori}, the algorithm must proceed from $r_{k-1}$ to $r_k$ by taking small increments $r_k(0)-r_{k-1}(0)=\epsilon_k$ and choosing $r_0(0)$ sufficiently small.  
\end{enumerate}
\end{remark}

To end this section, we prove that the sequence of approximations $\{\tilde\psi_k\}_{k=1}^\infty$ obtained from \eqref{eqn:the-approx} convergence uniformly to $\psi$.  
\begin{theorem}
Suppose that $\psi$ is defined on the cylinder $\Omega=\{(\theta,r):0\leq\theta\leq 2\pi, 0\leq r\leq\epsilon\}$ and let $\tilde\epsilon<\epsilon$ be chosen so that if $r:[0,2\pi]\rightarrow\real$ is a trajectory of \eqref{eqn:polar-red-r} with $r(0)\leq\tilde\epsilon$ then $r(\theta)<\epsilon$.  Let $\tilde\psi_k$ be defined as in \eqref{eqn:the-approx}, with step-size  $r_j(0)-r_{j-1}(0)=\epsilon_j:=\tfrac{1}{k+1}\tilde\epsilon$, for $j=0,1,\ldots,k$.  Then $\tilde\psi_k\rightarrow\psi$ uniformly in $\tilde\Omega=\{(\theta,r):0\leq\theta\leq 2\pi, 0\leq r\leq r_{k}(\theta)\}$.
\end{theorem}
\noindent\textbf{Proof.}  We first note that the existence of $\tilde\epsilon$ follows by Lyapunov stability of the exosystem.  By construction, $r_j(0)=\tfrac{j+1}{k+1}\tilde\epsilon$ for $j=0,1,\ldots,k$, and in particular $r_k(0)=\tilde\epsilon$, thereby rendering the domain $\tilde\Omega$ independent of $k$.  Also we note that, by shrinking $\tilde\epsilon$ if necessary, by Gronwall's lemma it follows that
\begin{equation}\label{eqn:est-gron}
|\tilde{r}(\theta)-\bar{r}(\theta)|\leq |\tilde{r}(0)-\bar r(0)|e^{K\theta},
\end{equation}
for all trajectories $\theta\mapsto \tilde r(\theta)$ and $\theta\mapsto \bar r(\theta)$ of \eqref{eqn:polar-red-r} such that $0\leq \tilde r(0),\bar r(0)\leq\tilde\epsilon$, where $K$ is a Lipschitz constant independent of $\theta$.

Now by definition, we have that $\Psi_j(\theta,\sigma) = \psi(\theta, r_{j-1}(\theta)+\sigma)$ and therefore
\[
\frac{\partial^{N+1}\Psi_j}{\partial\sigma^{N+1}}(\theta,\sigma) = \frac{\partial^{N+1}\psi}{\partial r^{N+1}}(\theta,r_{j-1}(\theta)+\sigma).
\]
Hence, by Taylor's theorem,
\begin{align*}
\Psi_j(\theta,\sigma) &= \Psi^N_{j}(\theta,\sigma) + \frac{1}{N!}\int_0^\sigma (\sigma-\tau)^N \frac{\partial^{N+1}\Psi_j}{\partial\sigma^{N+1}}(\theta,\tau)\,d\tau\\[2ex]
&= \Psi^N_{j}(\theta,\sigma) + \frac{1}{N!}\int_0^\sigma (\sigma-\tau)^N \frac{\partial^{N+1}\psi}{\partial r^{N+1}}(\theta,r_{j-1}(\theta)+\tau)\,d\tau.
\end{align*}
Let $C=\max_{(\theta,r)\in\Omega}\|\frac{\partial^{N+1}\psi}{\partial r^{N+1}}(\theta,r)\|$, and we note that $C$ exists by continuity of $\frac{\partial^{N+1}\psi}{\partial r^{N+1}}$ on the compact set $\Omega$.  We therefore have that   
\[
\|\Psi_j(\theta,\sigma) - \Psi^N_{j}(\theta,\sigma)\| \leq C\frac{\sigma^{N+1}}{(N+1)!}
\] 
provided $0\leq\sigma\leq r_j(\theta)-r_{j-1}(\theta)$ for $\theta\in[0,2 \pi]$, for all $j=0,1,\ldots,k$.  Now since $r_j(0)\leq \tilde\epsilon$ for $j=0,1,\ldots,k$, it follows by \eqref{eqn:est-gron} that
\[
|r_{j}(\theta)-r_{j-1}(\theta)| \leq |r_j(0)-r_{j-1}(0)|e^{K\theta} \leq \tfrac{1}{k+1}\tilde\epsilon e^{2\pi K}.
\]        

Therefore, given $(\theta,r)\in\tilde\Omega$, say that $r_{j-1}(\theta)\leq r < r_j(\theta)$ for some $j\in\{0,1,\ldots,k\}$, it follows that
\begin{align*}
\|\psi(\theta,r)-\tilde\psi_k(\theta,r)\|&= \|\Psi_j(\theta,r-r_{j-1}(\theta))-\psi_{j}(\theta,r-r_{j-1}(\theta))\|\\[2ex]
&= \|\Psi_j(\theta,r-r_{j-1}(\theta))-\Psi^N_{j}(\theta,r-r_{j-1}(\theta))\|\\[2ex]
&\leq \frac{C}{(N+1)!} (r-r_{j-1}(\theta))^{N+1}\\[2ex]
&\leq \frac{C}{(N+1)!}  \left(\tfrac{1}{k+1}\tilde\epsilon e^{2\pi K}\right)^{N+1}.
\end{align*}
Hence, $\|\psi(\theta,r)-\tilde\psi_k(\theta,r)\|\rightarrow 0$ as $k\rightarrow\infty$ uniformly in $\tilde\Omega$.  This completes the proof.\hfill$\blacksquare$

%==============================================
\section{Examples}\label{sec:examples}
In this section we present examples illustrating our method.
\begin{example}\label{exm:egg-func}
In this example we take a linear dynamical system of the form \eqref{eqn:gen-cm}, whose center manifold is easily computed, perform a nonlinear change of coordinates and arrive at a nonlinear system on which we apply our method.  The true solution for the nonlinear system is then readily available and we can compare the approximations produced by our method with the true soluiton.  Consider then the linear dynamical system
\begin{align*}
\dot x_1 &= x_2 + \tfrac{1}{2}w_1 + \tfrac{1}{2}w_2\\
\dot x_2 &= x_3 + \tfrac{1}{3}w_1 + \tfrac{2}{3}w_2\\
\dot x_3 &= -x_1 -\tfrac{1}{2}w_1 + \tfrac{1}{2}w_2\\
\dot w_1 &= -w_2\\
\dot w_2 &= w_1.
\end{align*}
The center manifold equation for this system in the unknown mapping $x=\phi(w)$ is
\[
\parder{\phi}{w}(w) Sw = C\phi(w) + Dw
\]
where $x=(x_1,x_2,x_3), w=(w_1,w_2)$, and
\[
S=\begin{pmatrix}0 & -1\\ 1 & 0\end{pmatrix}, C=\begin{pmatrix} 0 & 1 & 0\\0 & 0 & 1\\-1 & 0 & 0\end{pmatrix}, D=\begin{pmatrix} \tfrac{1}{2} & \tfrac{1}{2}\\[2ex] \tfrac{1}{3} & \tfrac{2}{3}\\[2ex] -\tfrac{1}{2} & \tfrac{1}{2}\end{pmatrix}.
\]
It is straightforward to verify that $\phi(w)=(-\tfrac{1}{3}w_1,-\tfrac{1}{2}w_1-\tfrac{1}{6}w_2,-\tfrac{1}{2}w_1-\tfrac{1}{6}w_2)$ is the unique solution to the center manifold equation for this system.  Consider the coordinate change $z=Z(x)=(-3x_1,9x_1-6x_2,-x_2+x_3+\rho(-3x_1,9x_1-6x_2))$, where $\rho:\real^2\rightarrow\real$ is a smooth function.  The system in $(z,w)$ coordinates takes the form of \eqref{eqn:gen-cm} with matrix $B$ having eigenvalues $-1, \tfrac{1}{2}\pm\sqrt{3}i$.  By direct substitution, the solution to the center manifold equation of the system in $(z,w)$ coordinates is $z=Z(\phi(w))=(w_1,w_2,\rho(w_1,w_2))$, i.e., it is the graph of the function $\rho$.  For purposes of illustration we take the egg carton shaped function $\rho(w_1,w_2)=\sin(w_1)\sin(w_2)$, whose graph is shown in Figure~\ref{fig:true-sol-egg}.  The patchy approximation $\tilde\psi_k(w_1,w_2)=(w_1,w_2,\tilde{\rho}(w_1,w_2))$ computed with our method with $k=10$ annular regions of thickness $\epsilon=0.5$ and of order $N=2$ is shown in Figure~\ref{fig:patchy-sol-deg7poly-egg}.  The error between the patchy approximation and the true solution is shown in Figure~\ref{fig:error-patchy-deg7poly-egg}.  Lastly, Figure~\ref{fig:error-poly-deg19-egg} shows that in order to get a similar error bound as with the patchy approximation, one needs to use a polynomial approximation of degree 19.  
\begin{figure}[thpb]
\centering
\includegraphics[width=12cm,keepaspectratio]{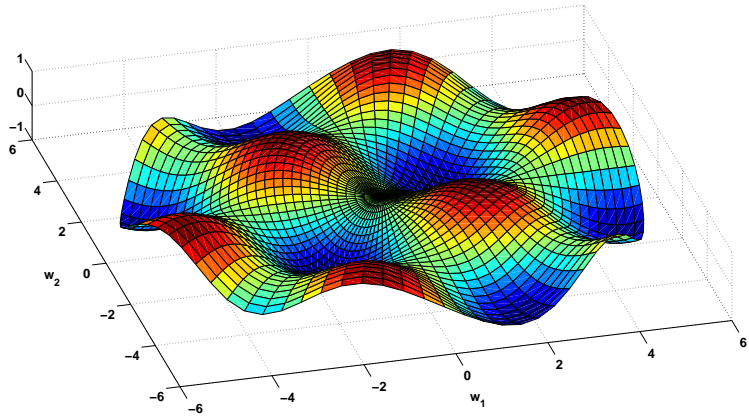}
\caption{True solution $\rho(w_1,w_2)=\sin(w_1)\sin(w_2)$.}
\label{fig:true-sol-egg}
\end{figure}

\begin{figure}[thpb]
\centering
\includegraphics[width=12cm,keepaspectratio]{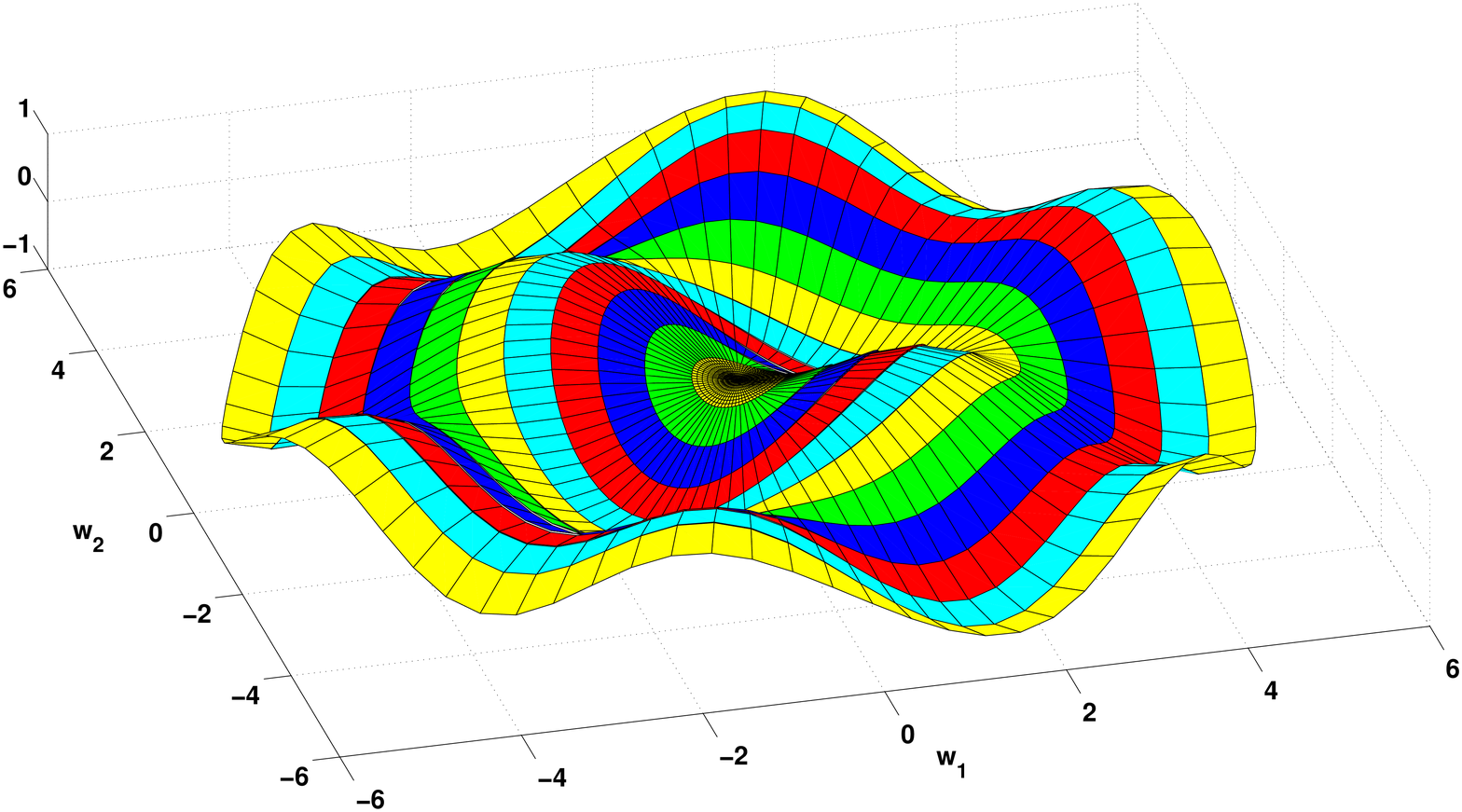}
\caption{Patchy solution $\tilde{\rho}(w_1,w_2)$ with $N=2$ and $k=10$.}
\label{fig:patchy-sol-deg7poly-egg}
\end{figure}

\begin{figure}[thpb]
\centering
\includegraphics[width=12cm,keepaspectratio]{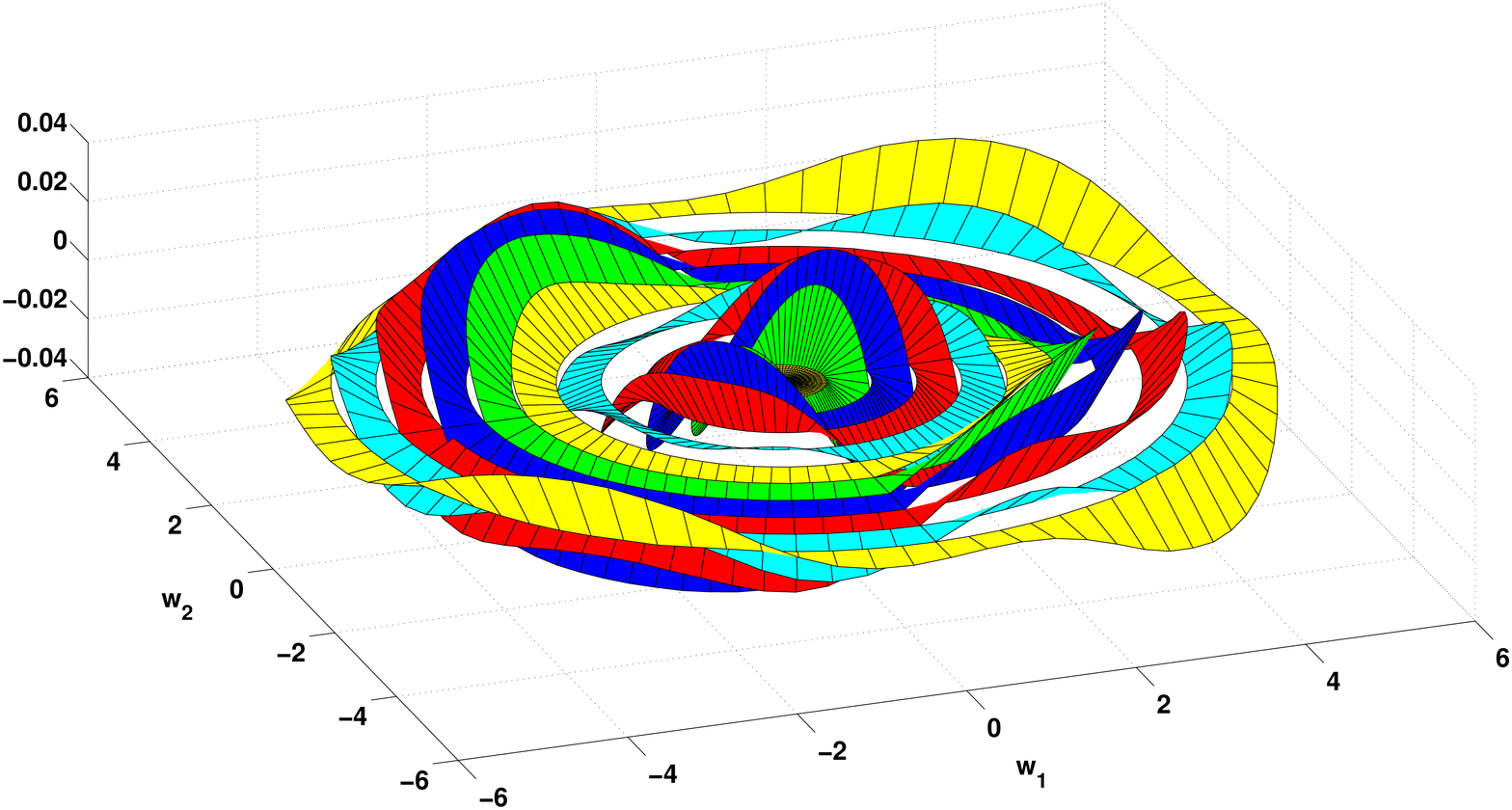}
\caption{Error $\rho(w_1,w_2)-\tilde{\rho}(w_1,w_2)$ with patchy solution with $N=2$ and $k=10$.}\label{fig:error-patchy-deg7poly-egg}
\end{figure}

\begin{figure}[thpb]
\centering
\includegraphics[width=12cm,keepaspectratio]{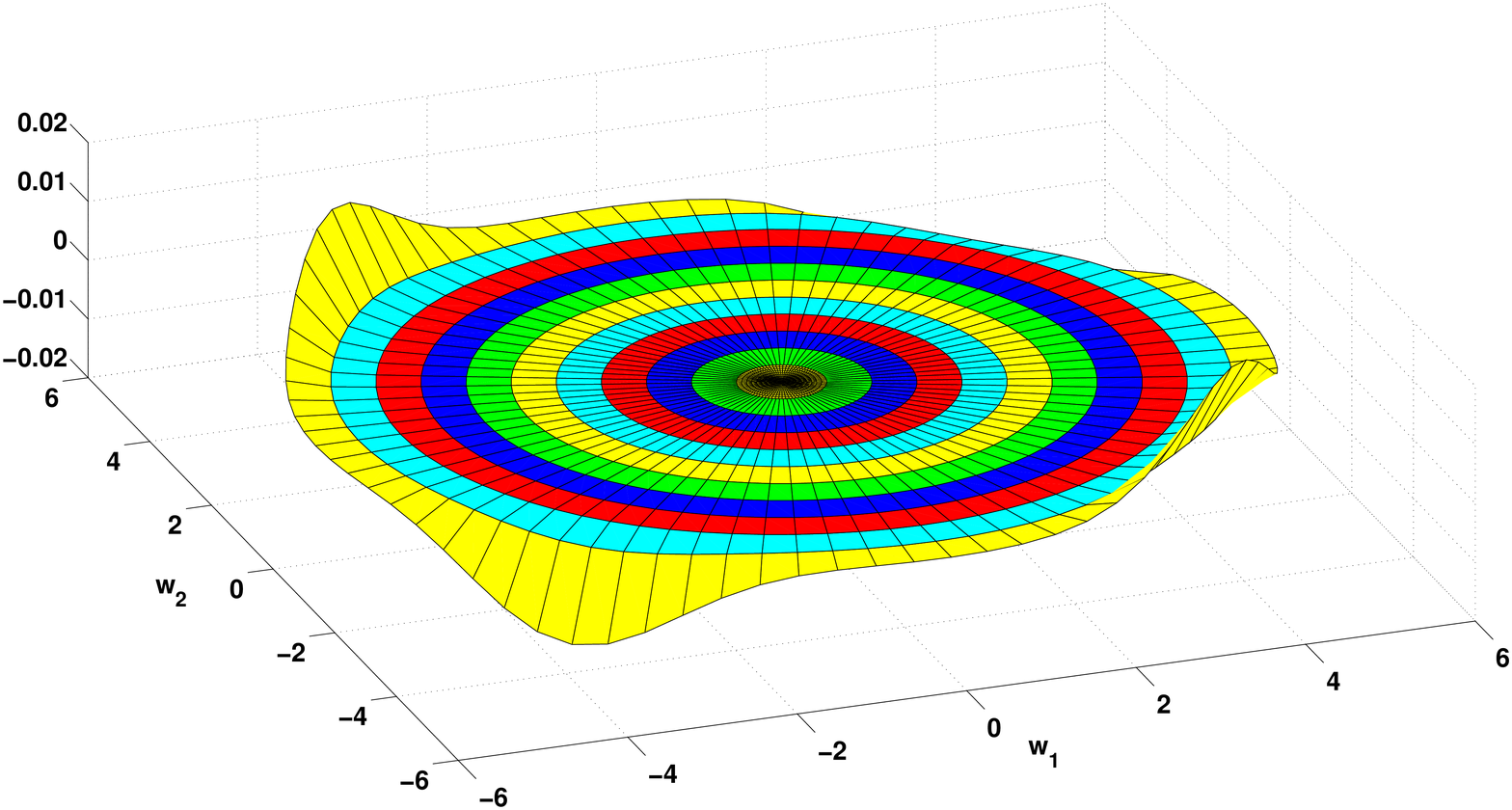}
\caption{Error $\rho(w_1,w_2)-\tilde{\rho}(w_1,w_2)$ with polynomial approximation of degree 19.}\label{fig:error-poly-deg19-egg}
\end{figure}

\end{example}

\begin{example}\label{exm:volcano-func}
This example illustrates the loss of stability when using polynomial approximations.  We proceed as in Example~\ref{exm:egg-func} but instead use the volcano type function $\rho(x,y)=\sin(x^2+y^2)e^{1-x^2-y^2}$, whose graph is shown in Figure~\ref{fig:true-sol-volcano}.  The patchy approximation is shown in Figure~\ref{fig:patchy-sol-volcano}  and the error in using the patchy approximation is shown in Figure~\ref{fig:error-patchy-volcano}.  The patchy approximation is constructed with $k=60$ annular regions of thickness $\epsilon=0.05$ and of order $N=1$, i.e., we only use a first order Taylor series in the radial direction.  For this example, polynomial approximations of orders up to $30$ where tested and it was verified that as one increases the order of the polynomial approximation the error in fact increases on the domain in consideration.

\begin{figure}[thpb]
\centering
\includegraphics[width=12cm,keepaspectratio]{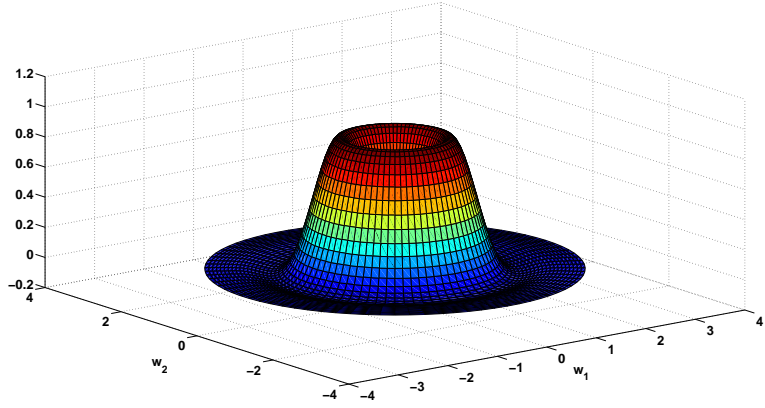}
\caption{True solution $\rho(x,y)=\sin(x^2+y^2)e^{1-x^2-y^2}$.}
\label{fig:true-sol-volcano}
\end{figure}

\begin{figure}[thpb]
\centering
\includegraphics[width=12cm,keepaspectratio]{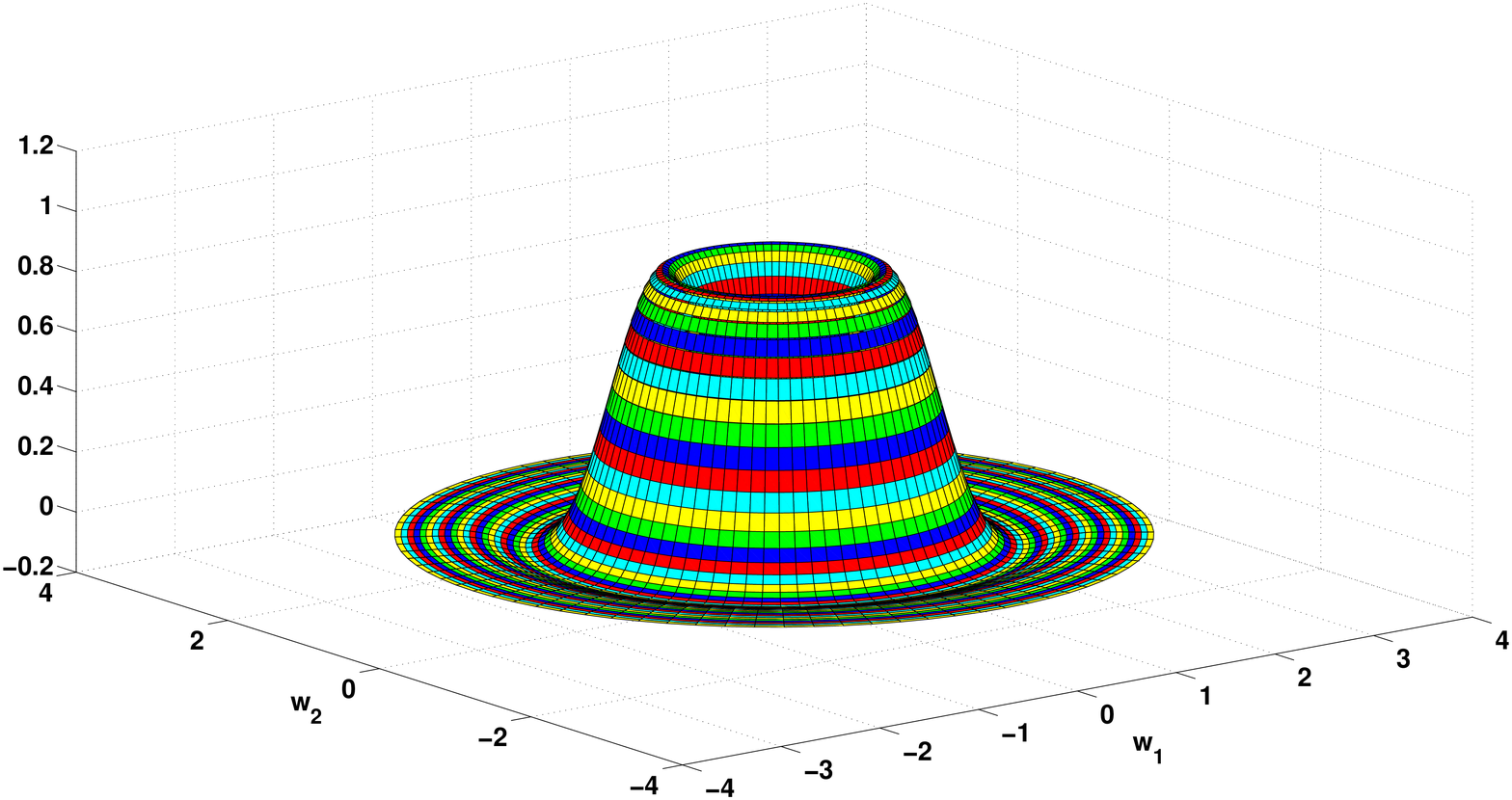}
\caption{Patchy solution $\tilde{\rho}(w_1,w_2)$ with $N=1$ and $k=60$.}
\label{fig:patchy-sol-volcano}
\end{figure}

\begin{figure}[thpb]
\centering
\includegraphics[width=12cm,keepaspectratio]{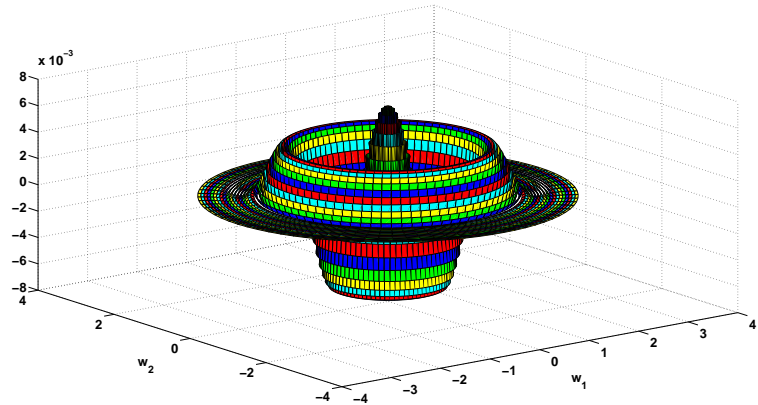}
\caption{Error $\rho(w_1,w_2)-\tilde{\rho}(w_1,w_2)$ with patchy solution with $N=1$ and $k=60$.}\label{fig:error-patchy-volcano}
\end{figure}

\end{example}

\begin{example}
Consider the inverted pendulum cart system from Example~\ref{exm:pendulum} with two-dimensional exosystem given by
\begin{equation}\label{eqn:duffing}
\dot w = s(w) = \begin{pmatrix} w_2\\ -w_1-aw_1^3\end{pmatrix}
\end{equation}
where $a>0$.  System~\eqref{eqn:duffing} is a special case of the unforced Duffing's oscillator with no damping \cite{GuckHol}.  The equilibrium $w=0$ of \eqref{eqn:duffing} is a center and representative periodic solutions encircling $w=0$ are shown in Figure~\ref{fig:duffing}.
\begin{figure}[thpb]
\centering
\includegraphics[width=12cm,keepaspectratio]{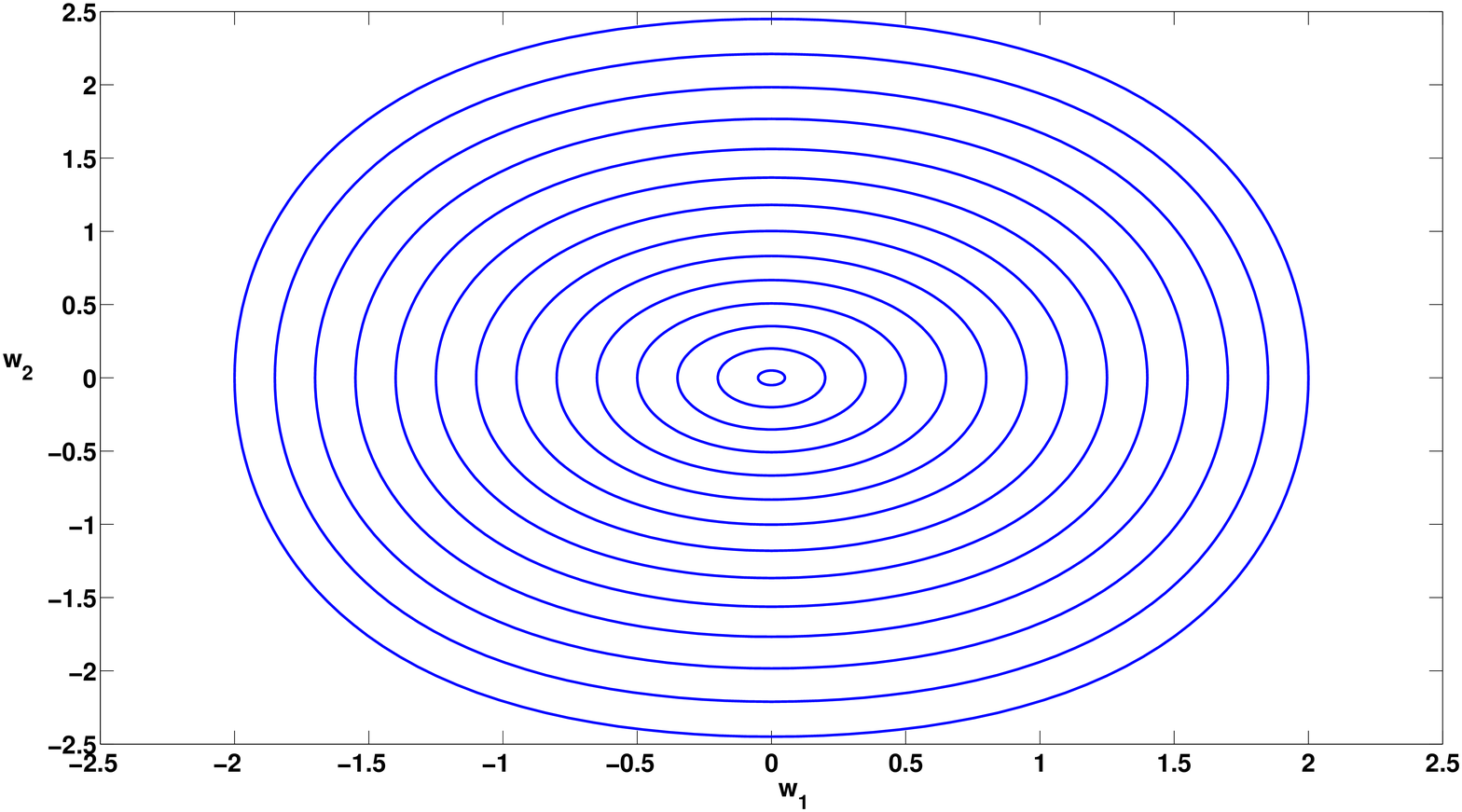}
\caption{Unforced Duffing's oscillator with no damping and $a=\tfrac{1}{4}$.}\label{fig:duffing}
\end{figure}

Let $p(w)=-w_1$.   As in Example~\ref{exm:pendulum}, let $(z,\xi)$ denote the standard normal coordinates, where $z=(x_3,x_4+\tfrac{x_2}{\ell}\cos(x_3))$ and $\xi=(h(x),L_fh(x))=(x_1,x_2)$.  Following the notation at the beginning of \S\ref{sec:ana-sol-fbi}, let $\varphi_1(w)=-p(w)=w_1$ and $\varphi_2(w)=-L_sp(w)=w_2$.  Then the differential equation \eqref{eqn:zero-exo} becomes
\begin{equation}\label{eqn:zero-pend-exo}
\begin{aligned}
\dot z_1 &= z_2 - \tfrac{1}{\ell}w_2\cos(z_1)\\
\dot z_2 &= \tfrac{g}{\ell}\sin(z_1) - \tfrac{1}{\ell}z_2w_2\sin(z_1) - \tfrac{1}{\ell^2}w_2^2\sin(z_1)\cos(z_1)\\
\dot w_1 &= w_2\\
\dot w_2 &= -w_1-aw_1^3
\end{aligned}
\end{equation}
A patchy approximation to the solution $\phi(w_1,w_2)=(\phi_1(w_1,w_2),\phi_2(w_1,w_2))$ of the center manifold equation of \eqref{eqn:zero-pend-exo} was computed and is illustrated in Figures~\ref{fig:phi1-pendulum}-\ref{fig:phi2-pendulum}.  The data $g=10$, $\ell=\tfrac{1}{3}$, and $a=\tfrac{1}{4}$ was used.  The solution was computed with $k=40$ and radial step at the initial angle $\theta=0$ was taken to be $\epsilon_i=0.05$, $i=1,2,\ldots,25$, and the order of the radial Taylor polynomials were chosen as $N=2$.

\begin{figure}[thpb]
\centering
\includegraphics[width=12cm,keepaspectratio]{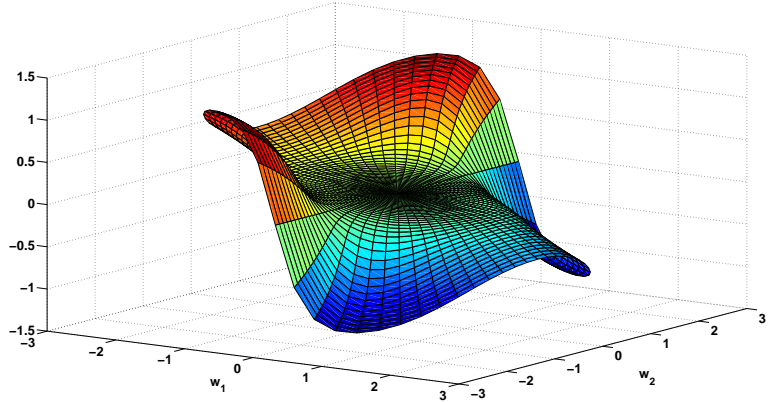}
\caption{Patchy approximation to $\phi_1(w_1,w_2)$.}\label{fig:phi1-pendulum}
\end{figure}

\begin{figure}[thpb]
\centering
\includegraphics[width=12cm,keepaspectratio]{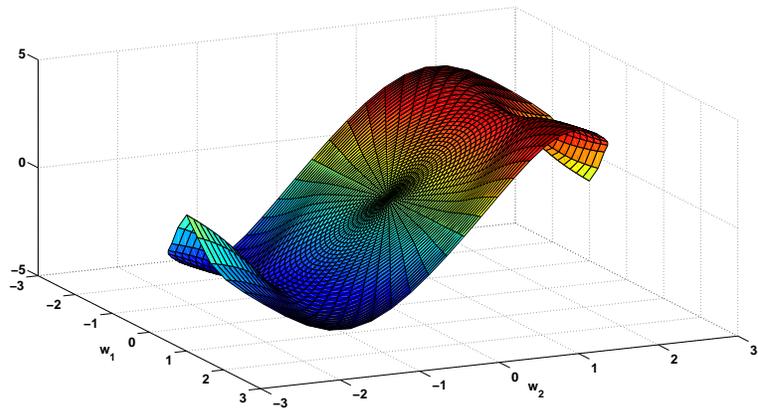}
\caption{Patchy approximation to $\phi_2(w_1,w_2)$.}\label{fig:phi2-pendulum}
\end{figure}

The computed patchy approximation to the center manifold PDE for \eqref{eqn:zero-pend-exo} is used in an output tracking controller of the form
\[
\alpha(z,w) = \kappa(w) + K((z,\xi)-\pi(w))
\]
where $\pi(w)=(\phi(w),\varphi(w))$ and $\kappa(w)=u_e(\pi(w),w)$, where the gain matrix $K$ is chosen as the solution to an LQR problem for the linearization of the inverted pendulum.  In the LQR problem, the matrices $Q=\textup{diag}(4,4,4,4)$ and $R=1$ were chosen.  A simulation is performed in which the pendulum is initialized at an angle of 15 degrees from the vertical and the cart is initialized at $-0.25$ from the origin.  The reference trajectory, $y_{\textup{ref}}(t)=w_1(t)$, is  chosen with initial condition $y_{\textup{ref}}(0)=1.2$.  The results of the simulation are shown in Figures~\ref{fig:tracking-pendulum}-\ref{fig:error-pendulum}.      

\begin{figure}[thpb]
\centering
\includegraphics[width=12cm,keepaspectratio]{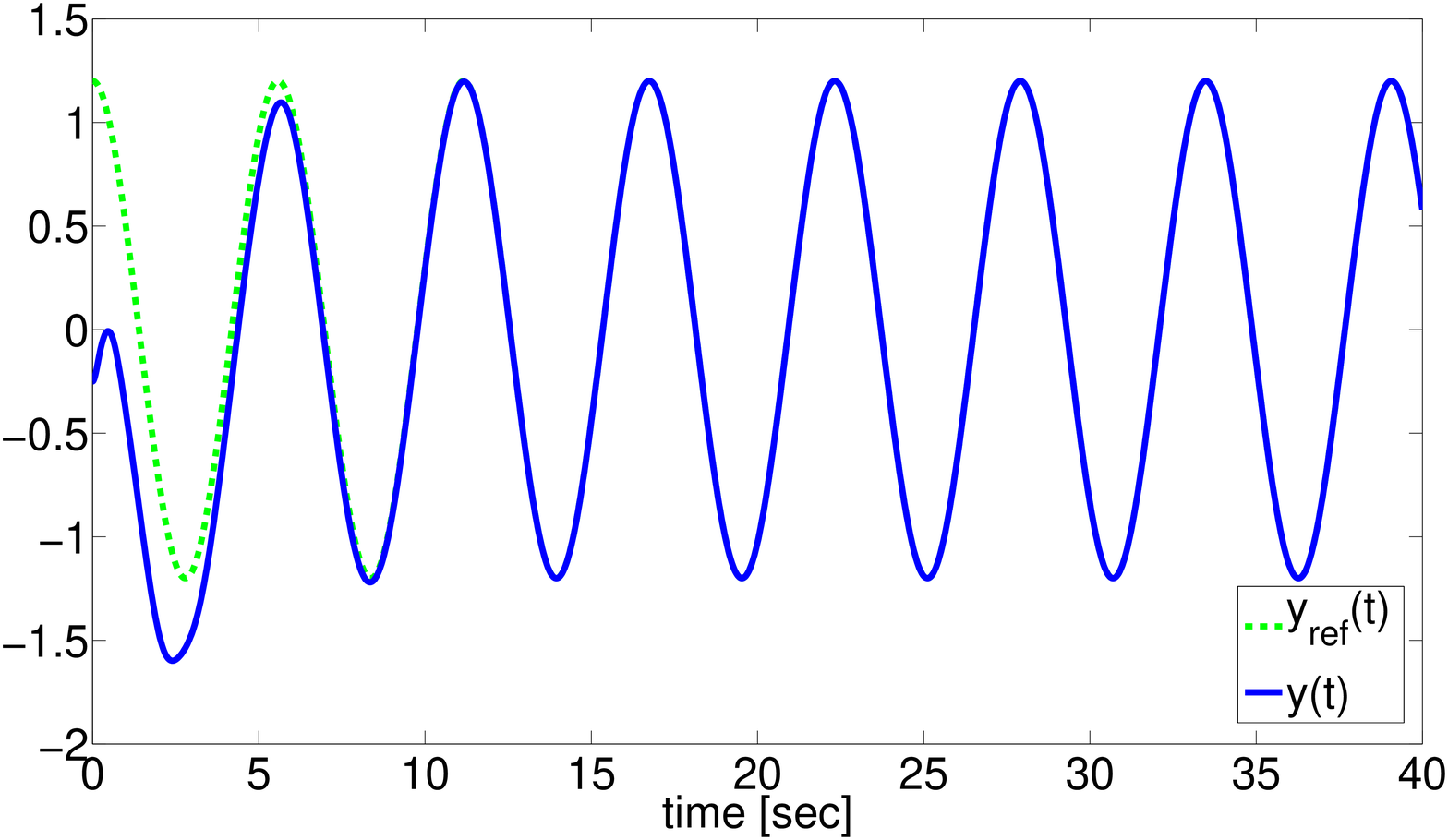}
\caption{Output $y(t)=x_1(t)$ and reference $y_{\textup{ref}}(t)=w_1(t)$.}
\label{fig:tracking-pendulum}
\end{figure}

\begin{figure}[thpb]
\centering
\includegraphics[width=12cm,keepaspectratio]{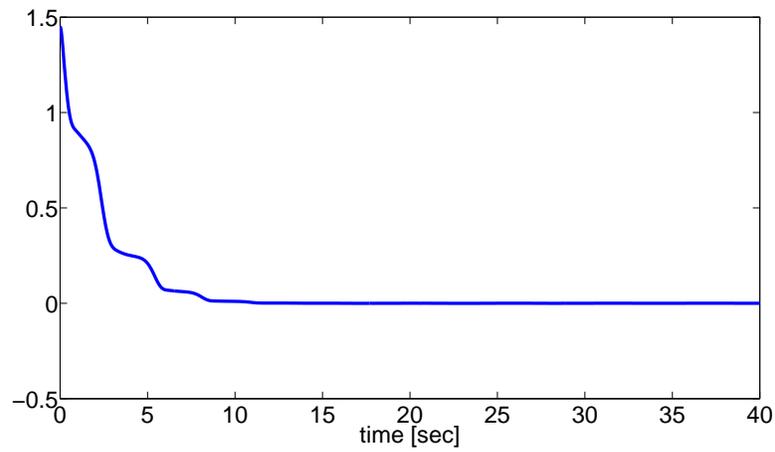}
\caption{Tracking error $e(t)=y(t)-y_{\textup{ref}}(t)$.}\label{fig:error-pendulum}
\end{figure}

\end{example}

%===================================
\section{Conclusion}
We have presented a method to compute solutions to the FBI equations of real analytic control-affine systems with two-dimensional exosystems.  Our technique is based on the patchy method in \cite{NavKr07} and on the results in \cite{Aulbach85} for uniqueness of solutions of two-dimensional real analytic center manifolds.  In comparison with direct Taylor polynomial approximations \cite{HuRu92,Kr92}, our method lessens the computational effort needed to produce approximate solutions by taking into account the periodic nature of a two-dimensional exosystem.  We proved that our method generates a sequence of approximations converging uniformly to the true solution.

\end{document}